\documentclass{article}

\usepackage{xcolor}
\usepackage{fourier}
\usepackage[normalem]{ulem}
\usepackage{amsmath, amsthm, dsfont}
\newtheorem{theorem}{Theorem}
\newtheorem{assumption}{Assumption}
\newtheorem{lemma}{Lemma}
\usepackage{authblk}

\usepackage[utf8]{inputenc}

\title{Consistency considerations for rare-event importance sampling}

\author[1]{J.\ Beh}
\author[2]{E.\ Flayac}
\author[3]{J.\ Morio}
\author[2]{F.\ Simatos}
\affil[1]{{\small Computing \& Mathematical Sciences Division,  Mohamed Bin Zayed University of Artificial Intelligence,  Abu Dhabi, UAE}}
\affil[2]{{\small Fédération ENAC ISAE-SUPAERO ONERA}}
\affil[3]{{\small ONERA/DTIS, Université de Toulouse, France}}

\def\d{\mathrm{d}}
\def\E{\mathds{E}}
\def\N{\mathds{N}}
\def\P{\mathds{P}}
\def\R{\mathds{R}}
\def\op{{\mathrm{op}}}
\def\det{\mathrm{det}}
\def\tr{\mathrm{tr}}
\def\Var{\mathds{V}\mathrm{ar}}
\newcommand{\Indicator}[1]{\mathds{1}\left(#1\right)}

\begin{document}

\maketitle

\abstract{We consider the problem of estimating a small probability $p = \P(X \in A)$ by importance sampling with auxiliary density $g$, where $X$ is drawn according to some density~$f$. We establish a simple necessary and sufficient condition for the consistency of the importance sampling estimator in triangular schemes where the dimension, the densities $f$ and $g$ and the set $A$ may evolve. It allows to refine the results of Chatterjee and Diaconis~\cite{Chatterjee18:0} by establishing the exact minimal sample size required for consistency. Moreover, it allows to study the regime $p \to 0$ which, despite being the relevant mathematical regime for rare-event estimations, has received very little attention so far. We perform an in-depth study of the hyperplane case $A = \{x: x^\top u \geq \alpha\}$ with Gaussian densities, showing that in this case, the optimal importance sampling density considerably improves consistency compared to the standard Monte Carlo estimator. However, we conclude the paper with moment considerations which suggest that, whenever the Monte Carlo and optimal importance sampling estimators both converge, convergence may occur faster for the Monte Carlo estimator.}

\section{Introduction}

\subsection{Triangular importance sampling schemes for rare-event estimation}

Let $Y \in \R^d$ be a random variable which is absolutely continuous with respect to Lebesgue measure with density $f$ (which we write $Y \sim f$) and $h: \R^d \to \R$ measurable so that $h(Y)$ is integrable. Importance sampling is a widely used simulation technique~\cite{Asmussen07:0, Elvira21:0} to estimate $\E[h(Y)]$ from samples from another distribution with density~$g$, called auxiliary distribution. The idea is simply to change measure to rewrite $\E[h(Y)] = \E[\ell(X) h(X)]$ with $X \sim g$ and $\ell = f/g$ the likelihood ratio between $f$ and $g$, and then estimate $\E[\ell(X) h(X)]$ by standard Monte Carlo with an i.i.d.\ sample drawn according to~$g$.

In this paper we are more specifically interested in the reliability context where one seeks to estimate a (small) probability $p = \P(Y \in A)$: this corresponds to the case where $h = \xi_A$, the indicator function of the set $A$, and the corresponding importance sampling estimator is then given by
\[ \hat p = \frac{1}{n} \sum_{i=1}^n \ell(X_i) \xi_A(X_i), \ X_i \sim g \ \text{i.i.d.}, \ \ell = \frac{f}{g}. \]
Note that for this estimator to be well defined with the usual convention $\frac{0}{0} = 0$, one needs $f\xi_A$ to be absolutely continuous with respect to $g$, i.e., $g(x) = 0 \Rightarrow f(x) \xi_A(x) = 0$. The conditional distribution $f|_A = f \xi_A / p$ is called the zero-variance auxiliary distribution since $\hat p = p$ when $g = f|_A$.  When $p$ is small, controlling the absolute error $\hat p - p$ is not necessarily insightful, and one rather turns to the relative estimator $\hat p / p$, leading to considering
\begin{equation} \label{eq:hat-p/p}
	\frac{\hat p}{p} = \frac{1}{np} \sum_{i=1}^n \ell(X_i) \xi_A(X_i) = \frac{1}{n} \sum_{i=1}^n \ell_A(X_i), \ X_i \sim g \ \text{i.i.d.}, \ \ell_A = \frac{\ell \xi_A}{p} = \frac{f \xi_A}{p g} = \frac{f|_A}{g}.
\end{equation}
The previous display shows that $\hat p / p$ can be seen as the empirical mean of the likelihood ratios~$\ell_A(X_i)$: thus, the key object of interest is the likelihood ratio $\ell_A = f|_A/g$ between the conditional distribution~$f|_A$ and the auxiliary distribution $g$, rather than the original likelihood ratio $\ell = f/g$ between $f$ and $g$.

Since $\E[\ell_A(X)] = 1$, the law of large numbers readily entails that, for fixed problem parameters $d$, $f$, $g$ and $A$, we have $\hat p / p \to 1$ almost surely as $n \to \infty$. However, there are different reasons for considering triangular schemes where, for instance, $A$ and/or $d$ are not fixed. In reliability for instance, one is interested in estimating small probabilities, motivating the study of $\hat p / p$ in the regime $p \to 0$. Moreover, there has recently been a surge of interest for high-dimensional importance sampling schemes~\cite{Agapiou17:0, Au03:0, Beh25-0, Uribe21:0}, corresponding to the regime $d \to \infty$. In such cases, we are in a triangular setting where the law of the summands $\ell_A(X_i)$ is not fixed and evolves together with the other problem parameters. The goal of this paper is to study this regime.

Both regimes have been studied in the literature, but mostly through moment-based criteria. In rare-event simulation, the regime $p \to 0$ is classically handled by requiring the relative second moment $\E[\ell_A(X)^2] = \E[\ell_A(Y_A)]$ to remain bounded (bounded relative error) or to satisfy $\log \E[\ell_A(Y_A)] = o(\log(1/p))$ (logarithmic efficiency); large deviations theory is then the main tool to design efficient auxiliary distributions, see for instance~\cite{Asmussen07:0, LEcuyer10:0, Sadowsky90:0} and~\cite{Glasserman97:0} for the refinements of this approach. In the high-dimensional regime $d \to \infty$, the degeneracy of importance sampling has been pointed out in reliability~\cite{Au03:0} and in filtering~\cite{Bengtsson08:0}, and the required sample size has been related to the second moment of the likelihood ratio~\cite{Agapiou17:0} or to the Kullback--Leibler divergence~\cite{Chatterjee18:0, SanzAlonso18:0}. In this paper we take a distributional viewpoint and focus on the consistency of $\hat p$, i.e., on the minimal sample size ensuring $\hat p / p \Rightarrow 1$.

\subsection{Triangular setting} \label{sub:triangular-scheme}

To formalize this, we consider in the rest of the paper a sequence of problems $((d_r, f_r, g_r, A_r, n_r), r \geq 1)$ where for each $r \geq 1$, $d_r \in \N \setminus \{0\}$ is the problem dimension, $A_r \subset \R^{d_r}$ is measurable and $f_r$ and $g_r$ are $d_r$-dimensional densities such that $f_r \xi_{A_r}$ is absolutely continuous with respect to $g_r$. Moreover, we define $p_r = \P(Y \in A_r)$, $Y \sim f_r$. To avoid cumbersome notation, the dependency on $r$ is omitted from the notation, but every limiting regime, supremum or infimum is taken with respect to $r$: for instance, $d \to \infty$ means that $d_r \to \infty$ as $r \to \infty$, $\inf p = \inf_r p_r$, etc. For simplicity, we will say in the sequel that $\hat p$ is consistent if $\hat p / p \Rightarrow 1$ with $\Rightarrow$ referring to convergence in distribution. Note that since the dimension $d$ may evolve (e.g., in the high-dimensional regime $d \to \infty$), the sequence of models is not naturally coupled: this is why it is more natural to only ask for convergence in distribution rather than almost surely.

\subsection{Notation}

In the rest of the paper, if $M$ is a square matrix, $\lVert M \rVert_\op$ denotes its operator norm, $\det(M)$ its determinant and $\tr(M)$ its trace. If $M$ is further symmetric, $\lambda_{\min}(M)$ and $\lambda_{\max}(M)$ denote its smallest and largest eigenvalue, respectively. Note that if $M$ is symmetric definite positive, then $\lVert M^{-1} \rVert_\op = 1/\lambda_{\min}(M)$ and $\lVert M \rVert_\op = \lambda_{\max}(M)$. Furthermore, with $f$, $g$ and $A$ as introduced in Section~\ref{sub:triangular-scheme}, $X$ and $X_i$ denote random variables with density $g$ and $Y$ and $Y_A$ denote random variables with density $f$ and $f|_A$, respectively. Moreover, for positive, deterministic sequences $u = (u_r)$ and $v = (v_r)$, we will write $u \ll v$ or $v \gg u$ to mean that $u/v \to 0$, or equivalently, $v/u \to +\infty$.

\section{Conditions for consistency, and relation with the Chatterjee--Diaconis approach}

\subsection{Closing the gap in the CD approach thanks to a simple necessary and sufficient condition}

As explained in the abstract of their groundbreaking paper~\cite{Chatterjee18:0}, Chatterjee and Diaconis have shown that, in a fairly general setting, a sample of size approximately $e^{D(f|_A, g)}$ is necessary and sufficient for accurate estimation by importance sampling, where $D(f|_A,g) = \E[\log (\ell_A(Y_A))] \in [0,\infty]$ is the Kullback--Leibler (KL) divergence between $f|_A$ and $g$, which exists provided that $f|_A$ is absolutely continuous with respect to $g$. In the sequel, we will refer for short to $n \gg e^{D(f|_A, g)}$ as the CD condition. Actually, although Chatterjee and Diaconis only mention the threshold $e^{D(f|_A, g)}$ in their abstract, their results are more precise. Indeed, \cite[Theorem~$1.3$]{Chatterjee18:0} states that for any $t \geq 0$ and $\delta \in (0,1)$,
\begin{equation} \label{eq:CD1}
	\E \left \lvert \frac{\hat p}{p} - 1 \right \rvert \leq e^{-t/4} + 2 \sqrt{\P \left( \log \ell_A(Y_A) > D(f|_A, g) + t/2 \right)} \ \text{ if } \ n = e^{D(f|_A, g) + t}
\end{equation}
and
\begin{equation} \label{eq:CD2}
	\P \left(\frac{\hat p}{p} \geq 1 - \delta \right) \leq e^{-t/2} + \frac{\P \left( \log \ell_A(Y_A) \leq D(f|_A, g) - t/2 \right)}{1-\delta} \ \text{ if } \ n = e^{D(f|_A, g) - t}.
\end{equation}
Roughly speaking, $t$ needs to grow faster than the fluctuations of the log-likelihood ratio around its mean for these two upper bounds to vanish. This suggests that if $t$ is larger than these fluctuations, then $\hat p$ will be consistent for $n = e^{D(f|_A, g) + t}$ and not consistent for $n = e^{D(f|_A, g) - t}$. To formalize this, we will sometimes consider the following assumption.

\begin{assumption} \label{ass:fluctuations}
	Let $L$ be a sequence of deterministic, positive numbers growing to $+\infty$. We say that the fluctuations of the log-likelihood ratio around its mean are of the order of $L$ if $D(f|_A, g) \gg L$ and if there exists a finite random variable $\bar L_\infty$ such that $\bar L \Rightarrow \bar L_\infty$, where
	\begin{equation} \label{eq:bar-L}
		\bar L = \frac{\log (\ell_A(Y_A)) - D(f|_A, g)}{L}.
	\end{equation}
\end{assumption}

When this assumption holds, then the CD bounds~\eqref{eq:CD1} and~\eqref{eq:CD2} readily imply the following result in terms of consistency.

\begin{lemma} \label{lemma:gap}
	Assume that the fluctuations of the log-likelihood ratio around its mean are of the order of $L$ in the sense of Assumption~\ref{ass:fluctuations}, and let $t \gg L$: then $\hat p$ is consistent if $n = e^{D(f|_A, g) + t}$ and not consistent if $n = e^{D(f|_A, g) - t}$.
\end{lemma}

Therefore, the CD approach leaves open the consistency of $\hat p$ in the range between $e^{D(f|_A, g) - L}$ and $e^{D(f|_A, g) + L}$: actually, there is a simple necessary and sufficient condition which allows to close this gap.

\begin{theorem}\label{thm:NSC}
	$\hat p / p \Rightarrow 1$ if and only if $\ell_A(Y_A) / n \Rightarrow 0$.
\end{theorem}

\begin{proof}
	Let $\chi^n_i = \ell_A(X_i) / n$, so that $\hat p / p = \sum_{i=1}^n \chi^n_i$ is a (triangular) sum of i.i.d.\ random variables. Since $\E(\chi^n_i) = 1/n$, we have $\chi_i^n \Rightarrow 0$ and so the triangular array $(\chi^n_i, 1 \leq i \leq n)$ is, in the language of~\cite[Definition VII.2.33]{Jacod03-0}, infinitesimal. Therefore, according to~\cite[Theorem VII.2.35]{Jacod03-0}, we have $\hat p / p \Rightarrow 1$ if and only if there exists a truncation function $h$, i.e., a bounded function equal to the identity function in a neighborhood of $0$, such that the following three conditions hold: 
	\begin{enumerate}
		\item[[$\beta$\!\!\!]] $n \E[h(\chi^n_1)] \to 1$;
		\item[[$\gamma$\!\!\!]] $n \left( \E[h(\chi^n_1)^2] - \{\E [h(\chi^n_1)]\}^2 \right) \to 0$;
		\item[[$\delta$\!\!\!]] $n \E[\phi(\chi^n_1)] \to 0$ for all bounded, continuous functions $\phi: \R \to \R$ which are $0$ around~$0$.
	\end{enumerate}
	In this case, [$\beta$] and [$\gamma$] hold for all truncation functions. Let us assume that $\hat p / p \Rightarrow 1$ and let us prove that $\ell_A(Y_A) / n \Rightarrow 0$. We use repeatedly the change of measure formula $\E [\varphi(Y_A)] = \E[ \ell_A(X) \varphi(X)]$ for any measurable, positive function $\varphi$. In particular, for any $\varepsilon > 0$ we have
	\[ \P(\ell_A(Y_A)/n \geq \varepsilon) = \E \left( \ell_A(X); \ell_A(X) / n \geq \varepsilon \right) = 1 - n \E[h(\chi_1^n)] \]
	with $h$ the truncation function $h(x) = x \Indicator{\lvert x \rvert < \varepsilon}$. Since $\hat p / p \Rightarrow 1$, we have $n \E[h(\chi^n_1)] \to 1$ by [$\beta$] giving $\ell_A(Y_A) / n \Rightarrow 0$ as desired.
	
	Let us now assume that $\ell_A(Y_A) / n \Rightarrow 0$, and let us prove that $\hat p / p \Rightarrow 1$, i.e., that [$\beta$], [$\gamma$] and [$\delta$] hold for $h(x) = x \Indicator{\lvert x \rvert < \varepsilon}$ for some $\varepsilon > 0$. For [$\beta$] this is the same argument as above, since $n \E[h(\chi^n_1)] = \P(\ell_A(Y_A)/n < \varepsilon)$. Let us now prove [$\gamma$]. Since $n \E h(\chi^n_1) \to 1$, we have $n \{\E[h(\chi^n_1)]\}^2 \to 0$ and so we only have to prove that $n \E[h(\chi^n_1)^2] \to 0$. Changing measure gives
	\[ n \E[h(\chi^n_1)^2] = \frac{1}{n} \E \left( \ell_A(X)^2 ; \ell_A(X) / n \leq \varepsilon \right) = \E[h(\ell_A(Y_A) / n)]. \]
	Since $\ell_A(Y_A) / n \Rightarrow 0$ by assumption and $h(\ell_A(Y_A) / n)$ is bounded, we get $\E [h(\ell_A(Y_A) / n)] \to 0$.
	Let us now prove [$\delta$]. Since $n \E(\phi(\chi^n_1)) = n \E(\phi(\chi^n_1) ; \chi^n_1 \geq \varepsilon) \leq n \sup \phi \P(\chi^n_1 \geq \varepsilon)$ it is enough to prove that $n \P(\chi^n_1 \geq \varepsilon) \to 0$. We have
	\[ n \P(\chi^n_1 \geq \varepsilon) = n \P(\ell_A(X) \geq \varepsilon n) = \E \left( \frac{1}{\ell_A(Y_A)/n} ; \frac{1}{\ell_A(Y_A)/n} \leq \frac{1}{\varepsilon} \right) = \E(Z; Z \leq 1/\varepsilon) \]
	with $Z = \frac{1}{\ell_A(Y_A) /n}$ which, by assumption, converges in distribution to $+\infty$. In particular, $Z \Indicator{Z \leq 1/\varepsilon} \Rightarrow 0$ and since $Z \Indicator{Z \leq 1/\varepsilon}$ is bounded by $1/\varepsilon$, the bounded convergence theorem implies that $\E(Z;Z\leq1/\varepsilon) \to 0$ which gives the desired result.
\end{proof}

Thanks to Theorem~\ref{thm:NSC}, the gap between the CD condition and the consistency of $\hat p$ can be understood: on the one hand, $e^{D(f|_A, g)} / n = \exp\{ \E[\log(\ell_A(Y_A)/n)] \}$ so that $n \gg e^{D(f|_A, g)}$ if and only if $\E[\log(\ell_A(Y_A)/n)] \to -\infty$, while on the other hand, Theorem~\ref{thm:NSC} states that $\hat p / p \Rightarrow 1$ if and only if $ \log(\ell_A(Y_A) / n) \Rightarrow -\infty$: thus, the difference between the CD condition and the consistency of $\hat p$ lies in the potential discrepancy between the convergence in distribution and in the mean of the normalized log-likelihood ratio $\log(\ell_A(Y_A)/n)$. Because $\log(\ell_A(Y_A))$ is not of constant sign, there is no systematic link between the behavior of its mean and of its distribution and the behavior of the mean is governed by the behavior of the positive and negative parts:
\begin{equation} \label{eq:positive-negative-parts}
	\E[\log(\ell_A(Y_A)/n)] = \E[\log_+(\ell_A(Y_A))] - \E[\log_-(\ell_A(Y_A))] - \log n
\end{equation}
with $\log_+ x = \Indicator{x \geq 1} \log x$ and $\log_- (x) =-  \Indicator{x \leq 1} \log x = \log_+(1/x)$. For the record we state the following lemma which stems from the convexity inequality $\log x \leq x-1$.

\begin{lemma}\label{lemma:negative-part}
	We have $\E[\log_-(\ell_A(Y_A))] \leq 1$.
\end{lemma}


This result shows that the potential discrepancy between the behavior in the mean and in distribution only owes to the behavior of the positive part of the log-likelihood ratio. In particular, this simple observation has the following consequence.

\begin{lemma} \label{lemma:log-n}
	If $\log n \gg D(f|_A, g)$, then $\hat p$ is consistent. In particular, if $\sup D(f|_A, g) < \infty$, then $\hat p$ is consistent for any $n \gg 1$.
\end{lemma}

\begin{proof}
	Assume that $\log n \gg D(f|_A, g)$, i.e., $\E[\log(\ell_A(Y_A))] / \log n \to 0$. By~\eqref{eq:positive-negative-parts} we have
	\[ \frac{\E[\log_+(\ell_A(Y_A))]}{\log n} =\frac{\E[\log(\ell_A(Y_A))]}{\log n} + \frac{\E[\log_-(\ell_A(Y_A))]}{\log n} \]
	and since $\E[\log_-(\ell_A(Y_A))] \leq 1$ this implies that $\E[\log_+(\ell_A(Y_A))] / \log n \to 0$. Since $\log_+(x) \geq 0$ this yields $\log_+(\ell_A(Y_A)) / \log n \Rightarrow 0$ which eventually implies that $\log_+(\ell_A(Y_A)) - \log n = \log n ( \log_+(\ell_A(Y_A))/\log n - 1) \Rightarrow -\infty$. Since $\log_-(\ell_A(Y_A))$ is tight (it is positive and its mean is bounded by $1$ by Lemma~\ref{lemma:negative-part}), this implies that $\log(\ell_A(Y_A)) - \log n \Rightarrow -\infty$, i.e., $\ell_A(Y_A) / n \Rightarrow 0$ and so we can conclude with Theorem~\ref{thm:NSC}.
\end{proof}

Of course, this result is far from being sharp as it settles the question in the regime $n = e^\omega$ with $\omega \gg D(f|_A,g)$ whereas the action takes place for $\omega \sim D(f|_A, g)$: however, the interest of the previous lemma lies in its simplicity, with no assumption whatsoever. Under Assumption~\ref{ass:fluctuations}, we can close the gap and get the precise minimal sample size required for the consistency of $\hat p$.

\begin{lemma} \label{lemma:D+L}
	Assume that the fluctuations of the log-likelihood ratio around its mean are of the order of $L$ in the sense of Assumption~\ref{ass:fluctuations} and let $n = e^{D(f|_A, g) + cL}$ for some $c \in \R$: then $\hat p$ is consistent if $\P(\bar L_\infty \geq c) = 0$ and not consistent if $\P(\bar L_\infty > c) > 0$. In particular, if $\bar L_\infty$ is unbounded above (i.e., $\P(\bar L_\infty \geq c) > 0$ for every $c \in \R$), then $\hat p$ is consistent if and only if $n = e^{D(f|_A, g) + t}$ for some $t \gg L$.
\end{lemma}

\begin{proof}
	Assume that Assumption~\ref{ass:fluctuations} holds and let $n = e^{D(f|_A, g) + cL}$ for some $c \in \R$: the question is whether $\ell_A(Y_A) / n \Rightarrow 0$ or not. Let $\varepsilon \in (0,1)$ and $c_- < c < c_+$ two continuity points of $\bar L_\infty$: since $\P\left( \ell_A(Y_A) \geq \varepsilon n \right) = \P\left( \bar L \geq (\log \varepsilon) / L + c \right)$ with $\bar L \Rightarrow \bar L_\infty$ and $L \to \infty$, we have for $r$ large enough $\P \left( \bar L \geq c_+ \right) \leq \P\left( \ell_A(Y_A) \geq \varepsilon n \right) \leq \P\left( \bar L \geq c_- \right)$ and so letting first $r \to \infty$ and then $c_- \uparrow c$ and $c_+ \downarrow c$ we obtain
	\[ \P \left( \bar L_\infty > c \right) \leq \liminf \P\left( \ell_A(Y_A) \geq \varepsilon n \right) \leq \limsup \P\left( \ell_A(Y_A) \geq \varepsilon n \right) \leq \P\left( \bar L_\infty \geq c \right) \]
	which proves the desired result.\\
	Assume now that $\bar L_\infty$ is unbounded above. Then the previous argument readily implies that $\hat p$ is consistent if $n = e^{D(f|_A, g) + t}$ for some $t \gg L$. If not, then by working along a subsequence we have $\sup (t / L) < \infty$ and along this subsequence, the previous arguments imply that $\P(\ell_A(Y_A) / n \geq \varepsilon) \not \to 0$, so that $\hat p$ is not consistent by Theorem~\ref{thm:NSC}.
\end{proof}

\subsection{Comparison to the standard Monte Carlo estimator}

Let $\hat p^{\text{MC}} = \frac{1}{n} \sum_{i=1}^n \xi_A(Y_i)$ with the $Y_i$'s i.i.d.\ $\sim f$ be the standard Monte Carlo estimator of $p$: a natural question is whether importance sampling improves the situation compared to the standard Monte Carlo approach. The ``improvement'' can take different forms depending on the performance metric considered: here we focus on consistency and we typically would like the consistency of $\hat p^{\text{MC}}$ to imply that of $\hat p$. When this occurs, i.e. $\hat p^{\text{MC}}$ consistent implies $\hat p$ consistent, we say that the importance sampling estimator $\hat p$ is more consistent than the standard Monte Carlo estimator~$\hat p^{\text{MC}}$. Theorem~\ref{thm:NSC} can help shed light on this question.

First, note that since $\ell_A(Y_A) = 1/p$ almost surely (a.s.) when $g = f$, Theorem~\ref{thm:NSC} readily entails that $\hat p^{\text{MC}} / p \Rightarrow 1$ if and only if $np \to \infty$. Next, since $\xi_A(Y_A) = 1$ a.s., we have
\[ \P(\ell_A(Y_A)/n \geq \varepsilon) = \P(\ell(Y_A) \geq \varepsilon np) = \P(g(Y_A) \leq f(Y_A) / (\varepsilon np)). \]
By Theorem~\ref{thm:NSC} we need this probability to vanish to have $\hat p$ consistent, and we can distinguish two cases:
\begin{itemize}
	\item if $\hat p^{\text{MC}}$ is consistent, i.e., $np \to \infty$, then $g(Y_A)$ may be smaller than $f(Y_A)$ while having $\hat p / p \Rightarrow 1$, but it should not be smaller than $f(Y_A) / (np)$: in this case we start from a consistent estimator ($\hat p^{\text{MC}}$) and $\hat p$ may remain consistent even if $g$ is worse than $f$ (i.e., puts less mass on $A$);
	\item but if $\hat p^{\text{MC}}$ is not consistent, i.e., $np \not \to \infty$, then $g(Y_A)$ must typically be larger than $f(Y_A)$ by a factor $1/(np)$ in order to ensure the consistency of $\hat p$: in this case, we start from a non-consistent estimator $\hat p^{\text{MC}}$ and provided $g$ puts more mass on $A$, importance sampling can yield a consistent estimator.
\end{itemize}

A natural way to capture whether $g$ puts more mass on $A$ than $f$ or not is by considering the probability $q := \P(X \in A)$ of the rare event set $A$ under the auxiliary density~$g$. We note that $nq \to \infty$ is a necessary condition for the consistency of $\hat p$ because the number of non-zero terms in~\eqref{eq:hat-p/p} follows a binomial distribution with parameter $(n,q)$ so that if $nq \to c \in [0,\infty)$, then $\P(\hat p = 0) = (1-q)^n \to e^{-c} > 0$ and we cannot have $\hat p / p \Rightarrow 1$. Based on this simple consideration, the following lemma shows that if $\hat p$ is more consistent than $\hat p^{\text{MC}}$, then the probability $q$ can, to the very least, not be negligible compared to $p$.

\begin{lemma}
	If $q \ll p$, then $\hat p$ is not more consistent than $\hat p^{\text{MC}}$.
\end{lemma}

\begin{proof}
	Assume that $q \ll p$ and let $n = 1/\sqrt{pq}$. Then $np = \sqrt{p/q} \to \infty$ so that $\hat p^{\text{MC}}$ is consistent, while $nq = \sqrt{q/p} \to 0$ so that $\hat p$ is not consistent.
\end{proof}

On the other hand, the CD condition suggests to consider as auxiliary density a density close to the conditional density $f|_A$ in the sense of the KL divergence: this is actually at the heart of the celebrated cross-entropy approach~\cite{Rubinstein04:0}, which consists in choosing $g$ minimizing $D(f|_A, g)$ among a given family of auxiliary distributions. The next lemma exhibits a simple relation between $q$ and $D(f|_A, g)$ which implies that if $g$ is closer to $f|_A$ than $f$ in the sense of the KL divergence, then necessarily $q \geq p$, i.e., $A$ has a larger probability under $g$ than under~$f$.

\begin{lemma} \label{lemma:q-D}
	We have $q \geq e^{-D(f|_A, g)}$. In particular:
	\begin{itemize}
		\item if $D(f|_A, g) \leq D(f|_A, f)$, then $q \geq p$;
		\item if $D(f|_A, g) \ll D(f|_A, f)$, then $\hat p$ is more consistent than $\hat p^{\text{MC}}$.
	\end{itemize}
\end{lemma}

\begin{proof}
	Let $g|_A = g \xi_A / q$ be the density $g$ conditioned on $A$. Since $f|_A$ is assumed to be absolutely continuous with respect to $g$, it is also absolutely continuous with respect to $g|_A$ and in particular, the two divergences $D(f|_A, g|_A)$ and $D(f|_A, g)$ are well-defined. Moreover, since $\xi_A(Y_A) = 1$ a.s., we have $g|_A(Y_A) = g(Y_A) / q$ which implies that
	\[ D(f|_A, g|_A) = \E \left[ \log \left(\frac{f|_A(Y_A)}{g|_A(Y_A)} \right) \right] = \E \left[ \log \left(\frac{q f|_A(Y_A)}{g(Y_A)} \right) \right] = \log q + D(f|_A, g). \]
	Since $D(f|_A, g|_A) \geq 0$ this implies the desired first inequality $q \geq e^{-D(f|_A, g)}$. The inequality $q \geq p$ then follows from the fact that $f|_A(Y_A) = f(Y_A) / p$ a.s., so that $D(f|_A, f) = -\log p$. In particular, if $D(f|_A, g) \leq D(f|_A, f)$, then $e^{-D(f|_A, f)} = p \leq e^{-D(f|_A, g)} \leq q$.
	
	To conclude, let us assume that $D(f|_A, g) \ll D(f|_A, f)$: if $\hat p^{\text{MC}}$ is consistent, then $n \gg 1/p = e^{D(f|_A, f)}$: therefore, $\log n \geq D(f|_A, f)$ for $r$ large enough, and so $\log n \gg D(f|_A, g)$ which implies the consistency of $\hat p$ by Lemma~\ref{lemma:log-n}.
\end{proof}

We have seen that $nq \to \infty$ (i.e., the mean number of non-zero terms in the sum~\eqref{eq:hat-p/p} defining $\hat p$) is a necessary condition for the consistency of $\hat p$. For the standard Monte Carlo estimator, this condition (which then reads $np \to \infty$) is also sufficient. A natural question is thus whether $nq \to \infty$ is also sufficient for consistency of $\hat p$: the next lemma, which is a direct consequence of Lemmas~\ref{lemma:D+L} and~\ref{lemma:q-D}, shows that this is not necessarily the case.

\begin{lemma}
	Assume that the fluctuations of the log-likelihood ratio around its mean are of the order of $L$ in the sense of Assumption~\ref{ass:fluctuations}, and that $\bar L_\infty$ is unbounded above. Then for $n = e^{D(f|_A, g) + t}$ with $1 \ll t \ll L$, $\hat p$ is not consistent while $nq \to \infty$.
\end{lemma}

\begin{proof}
	For $n = e^{D(f|_A, g) + t}$ with $1 \ll t \ll L$, $\hat p$ is not consistent by Lemma~\ref{lemma:D+L}, while $nq \geq e^t \to \infty$ by Lemma~\ref{lemma:q-D}.
\end{proof}

\section{Insight from the Gaussian setting} \label{sec:Gaussian}

Some of the above ideas can be further explored in the Gaussian case where explicit expressions are available. So in this section, we assume that $f$ is the density of the standard Gaussian distribution $N(0,I)$, and we also consider $g$ a Gaussian density, with mean $\mu$ and variance $\Sigma$.

In this setting, an important auxiliary distribution is the optimal (in the sense of the KL divergence) Gaussian density $g_A = \arg \min D(f|_A, g)$ where the minimum runs over all Gaussian densities. It is well-known that $g_A$ is the Gaussian density with mean $\mu_A := \E(Y_A)$ and covariance matrix $\Sigma_A := \Var(Y_A)$, see for instance~\cite{Beh25-0, Rubinstein04:0}. By construction, $g_A$ decreases the KL divergence with $f|_A$ and, based on the various results from the previous section, we expect it to improve the standard Monte Carlo estimator: we discuss some results quantifying this. 

\subsection{$g_A$ puts more mass on $A$}

In the Gaussian setting it is well-known that
\begin{equation} \label{eq:KL-gaussian-g_A}
	D(f|_A, g_A) = - \log p - \frac{1}{2} \mathcal{L}(\Sigma_A) - \frac{1}{2} \lVert \mu_A \rVert^2
\end{equation}
where $\mathcal{L}(M) = \tr(M-I) - \log \det(M) \geq 0$ is Stein's loss between $M$ and the identity matrix $I$, which satisfies $\mathcal{L}(M) = 0$ if and only if $M = I$, see for instance~\cite[Lemma~$3.5$]{Beh25-0}. Combined with Lemma~\ref{lemma:q-D}, this identity implies that, if $X \sim g_A$,
\[ \P(X \in A) \geq p e^{\frac{1}{2} \mathcal{L}(\Sigma_A) + \frac{1}{2} \lVert \mu_A \rVert^2}. \]
Thus, $g_A$ puts more mass on $A$ than $f$, which is coherent with the intuition that $g_A$ is a better auxiliary distribution than~$f$.

\subsection{Hyperplane case: consistency improvement}

By construction we have $D(f|_A, g_A) \leq D(f|_A, f)$ (which is also readily visible on~\eqref{eq:KL-gaussian-g_A} since $D(f|_A, f) = -\log p$) but the improvement brought by $g_A$ over $f$ is unclear. When $\inf p > 0$, then $D(f|_A, f) = - \log p$ is bounded and thus $D(f|_A, g_A)$ also: therefore in this case, both estimators have the same consistency condition, namely $n \gg 1$ by Lemma~\ref{lemma:log-n}.

The case $p \to 0$ is much more challenging and, although this regime is at the heart of rare-event simulation, the consistency of $\hat p$ for $g = g_A$ does not seem to have been considered in the literature so far, even in the one-dimensional case $d = 1$. In this case the condition for consistency of the standard Monte Carlo estimator is $n \gg 1/p$, and so a natural question is whether a suitable choice of $g$ may weaken this condition to $n \gg \kappa$ for some $\kappa \ll 1/p$. The KL-based arguments of the previous section suggest that this will be the case if $D(f|_A, g_A) \ll D(f|_A, f)$: in view of~\eqref{eq:KL-gaussian-g_A}, this requires a particular asymptotic behavior for $\mathcal{L}(\Sigma_A) + \lVert \mu_A \rVert^2$, namely that $\mathcal{L}(\Sigma_A) + \lVert \mu_A \rVert^2 \sim - 2 \log p$.

We do not know whether such an asymptotic relation holds for general sets $A$ and to get insight into this question, we consider the hyperplane case $A = \{x: x^\top u \geq \alpha \}$ with $\alpha \to \infty$ and $u$ a sequence of $d$-dimensional unit vectors ($\lVert u \rVert = 1$). This example may look quite specific, but it is justified for both theoretical and numerical reasons: from a theoretical perspective, it is a tractable case and from a numerical perspective, it is a standard test case considered in virtually any applied paper trying to assess the efficiency of some newly proposed importance sampling scheme, see for instance~\cite{Au03:0, Beh25-0, Tabandeh22:0}. In this case, it is easy to see that
\[ \mu_A = m u \ \text{ and } \ \Sigma_A = I + (\lambda - 1) uu^\top \ \text{ with } \ m = \E(N \mid N \geq \alpha) \ \text{ and } \ \lambda = \Var(N \mid N \geq \alpha) \]
where here and in the sequel, $N = Y^\top u$  follows a standard, one-dimensional normal distribution. In particular, it follows that
\begin{equation} \label{eq:formula-D}
	D(f|_A, g_A) = - \log p - \frac{1}{2} (\lambda - 1 - \log \lambda) - \frac{1}{2} m^2.
\end{equation}
Careful asymptotic expansions for $p$, $\lambda$ and $m$ imply that $D(f|_A, g_A)$ is bounded. More precisely, the following result holds.

\begin{lemma} \label{lemma:hyperplane}
	Consider the Gaussian setting and the hyperplane case $A = \{x: x^\top u \geq \alpha\}$ with $\alpha \to \infty$ and $u$ a sequence of unit vectors. Then we have $D(f|_A, g_A) \to \frac{1}{2} \log(2\pi/e)$. In particular, $\hat p$ for $g = g_A$ is consistent for any $n \gg 1$.
\end{lemma}

This result therefore shows that, as far as consistency is concerned, $g_A$ considerably improves the situation compared to the standard Monte Carlo estimator where the consistency condition is $n \gg 1/p$. Whether such an improvement holds for general sets $A$ represents in our view an interesting open question.

\begin{proof}[Proof of Lemma~\ref{lemma:hyperplane}]
	Using the Gaussian integration by parts formulas $\int_\alpha^\infty x^k e^{-x^2/2} \d x = \alpha^{k-1} e^{-\alpha^2/2} + (k-1) \int_\alpha^\infty x^{k-2} e^{-x^2/2} \d x$ for $k = 1,2$, we get
	\begin{equation} \label{eq:m-1}
		m = \E(N \mid N \geq \alpha) = \frac{1}{p} \int_\alpha^\infty \frac{x}{\sqrt{2\pi}} e^{-x^2/2} \d x = \frac{e^{-\alpha^2/2}}{p \sqrt{2\pi}}
	\end{equation}
	and
	\[ \E(N^2 \mid N \geq \alpha) = \frac{1}{p} \int_\alpha^\infty \frac{x^2}{\sqrt{2\pi}} e^{-x^2/2} \d x = \frac{1}{p \sqrt{2\pi}} \left( \alpha e^{-\alpha^2/2} + \int_\alpha^\infty e^{-x^2/2} \d x \right) = \alpha m + 1 \]
	so that
	\begin{equation} \label{eq:lambda}
		\lambda = \E(N^2 \mid N \geq \alpha) - m^2 = \alpha m + 1 - m^2.
	\end{equation}
	Thus, everything gets expressed in terms of $p$, whose asymptotic expansion is well-known. Actually, to get the leading term in $D(f|_A, g_A)$, we need to consider the following expansion:
	\begin{equation} \label{eq:p}
		p = \frac{e^{-\alpha^2/2}}{\sqrt{2 \pi}} \left( \frac{1}{\alpha} - \frac{1}{\alpha^3} + \frac{3}{\alpha^5} + o(\alpha^{-6}) \right).
	\end{equation}
	From this expansion we readily get $-\log p = \frac{\alpha^2}{2} + \log \alpha + \frac{1}{2} \log(2\pi) + o(1)$, while getting the expansions for the terms $\lambda$, $\log \lambda$ and $m^2$ appearing in~\eqref{eq:formula-D} is more involved. Starting from the asymptotic expansion~\eqref{eq:p} of $p$ and the expressions~\eqref{eq:m-1} and~\eqref{eq:lambda} of $m$ and $\lambda$ in terms of $p$, computation detailed below shows that
	\[ \lambda = o(1), \ \log \lambda = - 2 \log \alpha + o(1) \ \text{ and } \ m^2 = \alpha^2 + 2 + o(1). \]
	Plugging these expressions into~\eqref{eq:formula-D}, we obtain
	\[ D(f|_A, g_A) = \frac{\alpha^2}{2} + \log \alpha + \frac{1}{2} \log(2\pi) + \frac{1}{2} - \log \alpha - \frac{1}{2} \alpha^2 - 1 + o(1) = \frac{1}{2} \log(2\pi/e) + o(1) \]
	which is precisely the result announced. Let us now detail the computation for $m$, $\lambda$ and $\log \lambda$. Combining~\eqref{eq:m-1} and~\eqref{eq:p} gives $m = \alpha (1 - \alpha^{-2} + 3 \alpha^{-4} + o(\alpha^{-5}))^{-1}$	and so using $(1 - \varepsilon)^{-1} = 1 + \varepsilon + \varepsilon^2 + \varepsilon^3 + o(\varepsilon^3)$ with $\varepsilon = \alpha^{-2} - 3 \alpha^{-4} + o(\alpha^{-5})$ we obtain
	\begin{align*}
		m & = \alpha \left(1 + (\alpha^{-2} - 3 \alpha^{-4}) + (\alpha^{-2} - 3 \alpha^{-4})^2 + (\alpha^{-2} - 3 \alpha^{-4})^3 + o(\alpha^{-5}) \right)\\
		& = \alpha \left(1 + \alpha^{-2} - 3 \alpha^{-4} + \alpha^{-4} + o(\alpha^{-5}) \right)
	\end{align*}
	and so $m = \alpha + \alpha^{-1} - 2 \alpha^{-3} + o(\alpha^{-4})$. This allows to get the required expansion for $m^2$:
	\begin{align*}
		m^2 = \alpha^2 \left(1 + \alpha^{-2} - 2 \alpha^{-4} + o(\alpha^{-5}) \right)^2 = \alpha^2 \left(1 + 2 \alpha^{-2} - 4 \alpha^{-4} + \alpha^{-4} + o(\alpha^{-5}) \right)
	\end{align*}
	which gives $m^2 = \alpha^2 + 2 - 3 \alpha^{-2} + o(\alpha^{-3})$. Plugging in the previous expressions for $m$ and $m^2$ into the expression $\lambda = \alpha m + 1 - m^2$ of~\eqref{eq:lambda}, we get
	\[ \lambda = \alpha^2 + 1 - 2 \alpha^{-2} + 1 - \alpha^2 - 2 + 3 \alpha^{-2} + o(\alpha^{-3}) = \alpha^{-2} + o(\alpha^{-3}). \]
	This finally allows to get the expansion $\log \lambda = - 2 \log \alpha + o(1)$ which achieves the proof.
\end{proof}

In practice, $g_A$ is usually out of reach because the parameters $m$, $\lambda$ and $u$ are unknown. For this reason, modern algorithms resort to adaptive importance sampling~\cite{Beh25-1, Bugallo17:0, Papaioannou19:0, Rubinstein04:0, Uribe21:0} whereby these unknown parameters are first estimated in a preliminary step. This typically leads to considering auxiliary densities with mean $\mu$ and variance $\Sigma$ given by
\begin{equation} \label{eq:mu-Sigma}
	\mu = \gamma v \ \text{ and } \ \Sigma = I + (\kappa - 1) vv^\top \ \text{ with } \ \gamma \in \R, \ \kappa > 0 \ \text{ and } \ v \in \R^d \ \text{ with } \ \lVert v \rVert = 1.
\end{equation}
In recent algorithms, $\gamma$, $\kappa$ and $v$ are estimations of the true parameters $m$, $\lambda$ and $u$. Thus, it is interesting to study how the errors made in estimating these parameters propagates in the KL divergence: again, explicit formulas are available in the Gaussian setting and in the hyperplane case.

\begin{lemma}\label{lemma:g}
	Consider the Gaussian setting with $g$ the Gaussian density with mean $\mu$ and variance $\Sigma$ given by~\eqref{eq:mu-Sigma}, and consider also the hyperplane case $A = \{x: x^\top u \geq \alpha\}$ with $\alpha \to \infty$ and $u$ a sequence of unit vectors. Then we have
	\begin{equation} \label{eq:D}
		D(f|_A, g) = D(f|_A, g_A) + \frac{1}{2} \mathcal{L}(\Sigma^{-1} \Sigma_A) + \frac{1}{2} (\mu - \mu_A)^\top \Sigma^{-1} (\mu - \mu_A)
	\end{equation}
	with
	\begin{equation} \label{eq:L}
		\mathcal{L}(\Sigma^{-1} \Sigma_A) = (1-\langle u, v \rangle^2) (1 - \lambda)(\kappa^{-1} - 1) + \lambda \kappa^{-1} - 1 - \log \left( \lambda \kappa^{-1} \right)
	\end{equation}
	and
	\begin{equation} \label{eq:quadratic}
		(\mu - \mu_A)^\top \Sigma^{-1} (\mu - \mu_A) = \kappa^{-1} \left( \gamma - m \langle u, v \rangle \right)^2 + m^2 \left( 1 - \langle u, v \rangle^2 \right).
	\end{equation}
\end{lemma}

These expressions allow to recover the fact that $D(f|_A, g)$ will be small if $\Sigma$ is close to $\Sigma_A$ and $\mu$ to $\mu_A$, but they give precise quantifications to the extent to which we must have $\langle u, v \rangle \approx 1$, $\kappa \approx \lambda$ and $\gamma \approx m$. The value of the scalar product $\langle u, v \rangle$ seems particularly important as if we do not have $\langle u, v \rangle \to 1$, then we necessarily have $D(f|_A, g) \to \infty$.

\begin{proof} [Proof of Lemma~\ref{lemma:g}]
	The relation~\eqref{eq:D} is well-known, see for instance~\cite[Lemma~$3.5$]{Beh25-0}. Since $\Sigma_A = I + (\lambda - 1) uu^\top$ and $\Sigma = I + (\kappa - 1) vv^\top$, we get $\Sigma^{-1} = I + (\kappa^{-1} - 1) vv^\top$ and so
	\[ \Sigma^{-1} \Sigma_A = I + (\kappa^{-1} - 1) vv^\top + (\lambda - 1) uu^\top + (\kappa^{-1} - 1) (\lambda - 1) (v^\top u) v u^\top. \]
	In particular, $\tr(\Sigma^{-1} \Sigma_A - I) = \kappa^{-1} - 1 + \lambda - 1 + (\kappa^{-1} - 1)(\lambda - 1) \langle u, v \rangle^2$ and since $\det(\Sigma^{-1} \Sigma_A) = \lambda / \kappa$, we get
	\[ \mathcal{L}(\Sigma^{-1} \Sigma_A) = \kappa^{-1} - 1 + \lambda - 1 + (\kappa^{-1} - 1)(\lambda - 1) \langle u, v \rangle^2 - \log (\lambda / \kappa) \]
	which can be rearranged to get~\eqref{eq:L}. For the quadratic term, we have
	\[ (\mu - \mu_A)^\top \Sigma^{-1} (\mu - \mu_A) = (\gamma v - m u)^\top \left( I + (\kappa^{-1}-1) vv^\top \right) (\gamma v - mu) \]
	from which~\eqref{eq:quadratic} can be derived by standard algebraic manipulation.
\end{proof}

\section{Opening discussion: does $g_A$ improve the speed of convergence?} \label{sec:opening}

Although we have seen that $\hat p$ may significantly improve consistency compared to the standard Monte Carlo estimator $\hat p^{\text{MC}}$, it seems to us that an important, essentially open question is the following: when both estimators are consistent (e.g., if $\hat p$ is more consistent than $\hat p^{\text{MC}}$ and $np \to \infty$), does $\hat p$ converge faster than $\hat p^{\text{MC}}$? Since
\[ \zeta \left( \frac{\hat p}{p} - 1 \right) = \sum_{i=1}^n Z_i \ \text{ with } \ Z_i = \frac{\zeta}{n}(\ell_A(X_i) - 1), \]
determining the speed of convergence (i.e., the scaling parameters $\zeta$ that makes this sum converge) boils down to studying a triangular sum of i.i.d.\ random variables as in Theorem~\ref{thm:NSC}. However, the arguments in the proof of Theorem~\ref{thm:NSC} do not allow us to identify a simple condition for convergence.

In order to get insight into this question, we turn into looking at the moments $\E[\ell_A(X)^{1+\eta}] = \E[\ell_A(Y_A)^\eta]$ of the summands, which are meant to give indication on their tail behavior and thus, ultimately, on the behavior of their sum. In general, moments of $\ell_A(X)$ are difficult to study because of the indicator function of the set $A$ which, in principle, prohibits general explicit expressions. For these reasons, moments of $\ell_A(Y_A)$ are often controlled by neglecting the indicator function through the bound
\begin{equation} \label{eq:bound-moment}
	\E[\ell_A(Y_A)^\eta] \leq \frac{1}{p^{\eta+1}} \E[\ell(Y)^\eta].
\end{equation}
This bound is insightful when $\inf p > 0$, but is extremely rough when $p \to 0$: however, as mentioned above, what happens in this case is essentially unknown.\\

In the Gaussian setting, the next lemma shows that explicit expressions, which allow to capture the dependency on $A$, are available for small enough moments. Moreover, this result allows to precisely evaluate the gap in~\eqref{eq:bound-moment}. Here and in the sequel, we consider again the Gaussian setting of Section~\ref{sec:Gaussian} and we define $\bar \mu = \Sigma^{-1} \mu$ and $\Delta = \Sigma^{-1} - I${\color{red}.}

\begin{lemma} \label{lemma:moments-ell}
	Consider the Gaussian setting of Section~\ref{sec:Gaussian} and $\eta > 0$. If $\eta (\lVert \Sigma^{-1} \rVert_\op - 1) < 1$, then it holds that
	\begin{equation} \label{eq:moment-ell_A}
		\E [\ell_A(Y_A)^\eta] = \frac{1}{p^{\eta+1}} \E[\ell(Y)^\eta] \times \P \left( (I - \eta \Delta)^{-1/2} Y - \eta (I - \eta \Delta)^{-1} \bar \mu \in A \right)
	\end{equation}
	where
	\begin{equation} \label{eq:moment-ell}
		\E[\ell(Y)^\eta] = \exp \left( \eta D(f,g) + \frac{1}{2} \mathcal{L}(I - \eta \Delta) + \frac{\eta^2}{2} \bar \mu^\top (I - \eta \Delta)^{-1} \bar \mu \right).
	\end{equation}
	Moreover, $\E[\ell(Y)^\eta] = \infty$ if $\eta (\lVert \Sigma^{-1} \rVert_\op - 1) > 1$.
\end{lemma}

\begin{proof}
	We have $\E[\ell(Y)^\eta] = \int f^{\eta+1} / g^\eta$ and, starting from the expressions of $f$ and $g$, one can express the ratio $f^{\eta+1} / g^\eta$ as
	\[ \frac{f(y)^{\eta + 1}}{g(y)^\eta} = (2\pi)^{-d/2} [\det(\Sigma)]^{\eta/2} \exp \left( -\frac{1+\eta}{2} \left( y^\top (I - \delta \Sigma^{-1}) y + 2 \delta y^\top \Sigma^{-1} \mu - \delta \mu^\top \Sigma^{-1} \mu \right) \right) \]
	with $\delta = \eta / (1+\eta)$. This shows that $f^{\eta+1} / g^\eta$ will be integrable if $\lambda_{\min} (I - \delta \Sigma^{-1}) > 0$, and not integrable if $\lambda_{\min} (I - \delta \Sigma^{-1}) < 0$. Since $\lambda_{\min} (I - \delta \Sigma^{-1}) = \frac{1}{1+\eta} \left( 1 - \eta (\lVert \Sigma^{-1} \rVert_\op - 1) \right)$ this shows that $\ell(Y)^\eta$ is integrable if $\eta (\lVert \Sigma^{-1} \rVert_\op - 1) < 1$ and not integrable if $\eta (\lVert \Sigma^{-1} \rVert_\op - 1) > 1$. Let us now assume that $\eta (\lVert \Sigma^{-1} \rVert_\op - 1) < 1$ and proceed with establishing~\eqref{eq:moment-ell} and~\eqref{eq:moment-ell_A}. Let in the sequel $V = (I - \delta \Sigma^{-1})^{-1}$ and $z = \delta \bar \mu$: then we have
	\[ y^\top (I - \delta \Sigma^{-1}) y + 2 \delta y^\top \Sigma^{-1} \mu - \delta \mu^\top \Sigma^{-1} \mu = (y + Vz)^\top V^{-1} (y + Vz) - z^\top V z - \delta \mu^\top \Sigma^{-1} \mu \]
	so that
	\[ \frac{f(y)^{\eta + 1}}{g(y)^\eta} = [\det(\Sigma)]^{\eta/2} \exp \left( \frac{1+\eta}{2} \left( z^\top V z + \delta \mu^\top \Sigma^{-1} \mu \right) \right) [\det(V/(1+\eta))]^{1/2} h(y) \]
	with $h$ the Gaussian density with variance $V / (1+\eta) = (I - \eta \Delta)^{-1}$ and mean $-Vz = -\eta (I - \eta \Delta)^{-1} \bar \mu$. Integrating this identity implies that
	\[ \E[\ell(Y)^\eta] = [\det(\Sigma)]^{\eta/2} \exp \left( \frac{1+\eta}{2} \left( z^\top V z + \delta \mu^\top \Sigma^{-1} \mu \right) \right) [\det(V/(1+\eta))]^{1/2} \]
	as well as~\eqref{eq:moment-ell_A} since $\E[\ell_A(Y_A)^\eta] = \frac{1}{p^{\eta+1}} \int \xi_A \frac{f^{\eta+1}}{g^\eta}$. Therefore, in order to conclude the proof it remains to prove~\eqref{eq:moment-ell}, i.e., in view of the above, that
	\begin{multline*}
		\frac{\eta}{2} \log \det(\Sigma) - \frac{1}{2} \log \det(I - \eta \Delta) + \frac{1+\eta}{2} \left( z^\top V z + \delta \mu^\top \Sigma^{-1} \mu \right)\\
		= \eta D(f,g) + \frac{1}{2} \mathcal{L}(I - \eta \Delta) + \frac{\eta^2}{2} \bar \mu^\top (I - \eta \Delta)^{-1} \bar \mu.
	\end{multline*}
	This relation is easily established recalling that $D(f,g) = \frac{1}{2} \mathcal{L}(\Sigma^{-1}) + \frac{1}{2} \mu^\top \Sigma^{-1} \mu$ and $\mathcal{L}(M) = \tr(M-I) - \log \det(M)$.
\end{proof}

Of course, depending on the set $A$, we may have $\E[\ell_A(Y_A)^\eta] < \infty$ even if $\E[\ell(Y)^\eta] = \infty$: for instance, consider $\eta (\lVert \Sigma^{-1} \rVert_\op - 1) > 1$, so that $\E[\ell(Y)^\eta] = \infty$ by the previous lemma, but $A = \{ x: \lVert x \rVert \leq 1 \}$ in which case $\ell_A$ is bounded and so $\ell_A(Y_A)$ has finite moments of all order. An inspection of the proof of Lemma~\ref{lemma:moments-ell} reveals that what happens is that $\E[\ell_A(Y_A)^\eta] \propto \int \xi_A(y) e^{-Q_\eta(y)} \d y$ with $Q_\eta$ a symmetric quadratic form. The condition $\eta (\lVert \Sigma^{-1} \rVert_\op - 1) < 1$ ensures that $Q_\eta$ is positive, which makes $e^{-Q_\eta}$ (and thus $\xi_A e^{-Q_\eta}$) integrable. If $\eta (\lVert \Sigma^{-1} \rVert_\op - 1) > 1$ then $Q_\eta$ has negative eigenvalues, which makes $e^{-Q_\eta}$ not integrable since it diverges in the directions of the eigenvectors associated to the negative eigenvalues. However, this diverging behavior could be compensated by the indicator function $\xi_A$: for instance, one can easily build $Q = -uu^\top + Q_+$ with $Q_+$ positive, such that $e^{-Q}$ is not integrable while $\xi_A e^{-Q}$ with $A = \{x: \lvert u^\top x \rvert \leq 1\}$ is integrable.

However, the hyperplane case suggests that such a ``cooperative'' interplay between $A$ and $g$ may not necessarily be the case: indeed, in this case the indicator function does not change the tail behavior and so the condition for the existence of moments of $\ell_A(Y_A)$ is the same as for $\ell(Y)$.

\begin{lemma} \label{lemma:moments-ell_A}
	Consider the Gaussian setting of Section~\ref{sec:Gaussian} with $g = g_A$, and the hyperplane case $A = \{x: x^\top u \geq \alpha\}$ with $u$ a unit vector and $\alpha \in \R$. For $\eta > 0$, we have $\E[\ell_A(Y_A)^\eta] < \infty$ if $\eta (1/\lambda - 1) < 1$ and $\E[\ell_A(Y_A)^\eta] = \infty$ if $\eta (1/\lambda - 1) > 1$.
\end{lemma}

This result is puzzling as, in the hyperplane case, we have $\lambda \to 0$ when $\alpha \to \infty$: therefore, although $g = g_A$ is a very good choice for consistency (cf.\ Lemma~\ref{lemma:hyperplane}), it seems extremely bad in terms of tail behavior as $\ell_A(Y_A)$ for $g = g_A$ has very few moments (up to $\lambda \sim \alpha^{-2}$) while the Monte Carlo estimator has all moments, since $\ell_A(Y_A) = 1/p$ for $g = f$. Based on these considerations, and given that the hyperplane is usually a favorable test case, we conjecture that, although $g_A$ dramatically improves consistency, in the regime $n p \to \infty$ where both $\hat p^{\text{MC}}$ and $\hat p$ with $g = g_A$ converge, then convergence may occur faster for $\hat p^{\text{MC}}$ than for $\hat p$.

\providecommand{\bysame}{\leavevmode\hbox to3em{\hrulefill}\thinspace}
\providecommand{\MR}{\relax\ifhmode\unskip\space\fi MR }
\providecommand{\MRhref}[2]{%
  \href{http://www.ams.org/mathscinet-getitem?mr=#1}{#2}
}
\providecommand{\href}[2]{#2}


\begin{thebibliography}{10}

\bibitem{Agapiou17:0}
Sergios Agapiou, Omiros Papaspiliopoulos, Daniel Sanz-Alonso, and Andrew~M.
  Stuart, \emph{Importance sampling: intrinsic dimension and computational
  cost}, Statist. Sci. \textbf{32} (2017), no.~3, 405--431.

\bibitem{Asmussen07:0}
S{\o}ren Asmussen and Peter~W. Glynn, \emph{Stochastic simulation: algorithms
  and analysis}, Stochastic Modelling and Applied Probability, vol.~57,
  Springer, New York, 2007.

\bibitem{Au03:0}
Siu-Kui Au and James~L. Beck, \emph{Importance sampling in high dimensions},
  Struct. Saf. \textbf{25} (2003), no.~2, 139--163.

\bibitem{Beh25-1}
Jason Beh, Jerome Morio, and Florian Simatos, \emph{Phase transition for
  conditional covariance matrices estimated by importance sampling, and
  implications for cross-entropy schemes in high dimension}, 2025, Preprint,
  arXiv:2511.11351.

\bibitem{Beh25-0}
Jason Beh, Yonatan Shadmi, and Florian Simatos, \emph{{Insight from the
  Kullback--Leibler divergence into adaptive importance sampling schemes for
  rare event analysis in high dimension}}, Ann. Appl. Probab. \textbf{35}
  (2025), no.~2, 1083--1124.

\bibitem{Bengtsson08:0}
Thomas Bengtsson, Peter Bickel, and Bo~Li, \emph{Curse-of-dimensionality
  revisited: collapse of the particle filter in very large scale systems},
  Probability and statistics: essays in honor of {D}avid {A}. {F}reedman, Inst.
  Math. Stat. Collect., vol.~2, Inst. Math. Statist., Beachwood, OH, 2008,
  pp.~316--334.

\bibitem{Bugallo17:0}
M{\'o}nica~F. Bugallo, V{\'\i}ctor Elvira, Luca Martino, David Luengo,
  Joaqu{\'\i}n M{\'\i}guez, and Petar~M. Djuri{\'c}, \emph{Adaptive importance
  sampling: the past, the present, and the future}, IEEE Signal Process. Mag.
  \textbf{34} (2017), no.~4, 60--79.

\bibitem{Chatterjee18:0}
Sourav Chatterjee and Persi Diaconis, \emph{The sample size required in
  importance sampling}, Ann. Appl. Probab. \textbf{28} (2018), no.~2,
  1099--1135.

\bibitem{Elvira21:0}
V{\'\i}ctor Elvira and Luca Martino, \emph{Advances in importance sampling},
  Wiley StatsRef: Statistics Reference Online, Wiley, 2021.

\bibitem{Glasserman97:0}
Paul Glasserman and Yashan Wang, \emph{Counterexamples in importance sampling
  for large deviations probabilities}, Ann. Appl. Probab. \textbf{7} (1997),
  no.~3, 731--746.

\bibitem{Jacod03-0}
Jean Jacod and Albert~N. Shiryaev, \emph{Limit theorems for stochastic
  processes}, second ed., Grundlehren der Mathematischen Wissenschaften
  [Fundamental Principles of Mathematical Sciences], vol. 288, Springer-Verlag,
  Berlin, 2003. \MR{MR1943877 (2003j:60001)}

\bibitem{LEcuyer10:0}
Pierre L'Ecuyer, Jose~H. Blanchet, Bruno Tuffin, and Peter~W. Glynn,
  \emph{Asymptotic robustness of estimators in rare-event simulation}, ACM
  Trans. Model. Comput. Simul. \textbf{20} (2010), no.~1, Article 6.

\bibitem{Papaioannou19:0}
Iason Papaioannou, Sebastian Geyer, and Daniel Straub, \emph{Improved cross
  entropy-based importance sampling with a flexible mixture model}, Reliab.
  Eng. Syst. Saf. \textbf{191} (2019), 106564.

\bibitem{Rubinstein04:0}
Reuven~Y. Rubinstein and Dirk~P. Kroese, \emph{The cross-entropy method: a
  unified approach to combinatorial optimization, {M}onte-{C}arlo simulation
  and machine learning}, Information Science and Statistics, Springer-Verlag,
  New York, 2004.

\bibitem{Sadowsky90:0}
John~S. Sadowsky and James~A. Bucklew, \emph{On large deviations theory and
  asymptotically efficient {M}onte {C}arlo estimation}, IEEE Trans. Inform.
  Theory \textbf{36} (1990), no.~3, 579--588.

\bibitem{SanzAlonso18:0}
Daniel Sanz-Alonso, \emph{Importance sampling and necessary sample size: an
  information theory approach}, SIAM/ASA J. Uncertain. Quantif. \textbf{6}
  (2018), no.~2, 867--879.

\bibitem{Tabandeh22:0}
Armin Tabandeh, Gaofeng Jia, and Paolo Gardoni, \emph{A review and assessment
  of importance sampling methods for reliability analysis}, Struct. Saf.
  \textbf{97} (2022), 102216.

\bibitem{Uribe21:0}
Felipe Uribe, Iason Papaioannou, Youssef~M. Marzouk, and Daniel Straub,
  \emph{Cross-entropy-based importance sampling with failure-informed dimension
  reduction for rare event simulation}, SIAM/ASA J. Uncertain. Quantif.
  \textbf{9} (2021), no.~2, 818--847.

\end{thebibliography}
\end{document}